\documentclass[a4paper,10pt]{article}
\def \Z {{\mathbf {Z}}}

\def \N {{\mathbf {N}}}

\def\eps{\varepsilon}
\title{ Lebesgue Component in Spectrum \\ for Tensor Product  of Generic  Multiplies}
\author{V.V. Ryzhikov}
\date{vryzh@mail.ru}

\begin{document}
\Large
\maketitle
{ \large We present a simple method to construct two  so-called generic measure-preserving transformations with  Lebesgue component in  spectrum of their tensor product.
We show how to get two rigid Gaussian systems and  two rigid  Poisson suspensions  with similar spectral properties.
 "Generic"    means here,  for instance, that a system $T$ is rigid, rank-one,   has all  polynomials   in its weak close.  Generic properties imply  simple spectrum  of Gaussian suspension $e^T$ ($\exp (T)$). 
({\large
 We call  a sums $\sum_{z\in \Z} a_z T^z$  "polynomial"  presuming  $a_z\geq 0$,  and   $\sum_z a_z=1$ in a case of  a probability space. For infinite spaces one assumes that $a_z\geq 0$  satisfy the condition  $\sum_z a_z\leq 1$.})}

To explain our simple idea let's  construct for now  ergodic  transformations  $S$, $T$  of an infinite measure space  such that

(i) $S$ and $T$ are rigid,\ \ \ 
(ii)  $S\otimes T$  has Lebesgue spectrum,

(iii) Poisson suspensions $S_\ast$, $T_\ast$  and Gaussian automorphisms
  $e^S, e^T$  have simple (hence, singular) spectrum. 

\bf Constructions.  \rm We  build  $T$ and $S$   defining  step by step  certain integer time-intervals  
\vspace{3mm}

\ \ \ \ \ \ \ \ \ \ \ \ \ $I_1$ \ \  $\tilde I_1$  \ \  $I_2$  \ \ $\tilde I_2$   \ \    $I_3$  \ \  \dots   $I_n$   \ \  $\tilde I_n$   \ \  $I_{n+1}\dots$,  \ \ \

\ \ \ \ \ \ \ \ \ \ \ \ \ $\tilde J_1$ \    $J_1$   \ \  $\tilde J_2$   \ \    $J_2$  \ \   $\tilde J_3$   \  \dots $\tilde  J_n$  \ \  $J_n$  \ \  $\tilde  J_{n+1}\dots, $\ \ \  

satisfying  
$$ \tilde I_n\subset J_n, \ \ \  \tilde J_n\subset I_n, \ \ \ \ \ \
\bigcup_n  (I_n\cup J_n)=\N.  $$  
\\ \\ \\ \\ \\
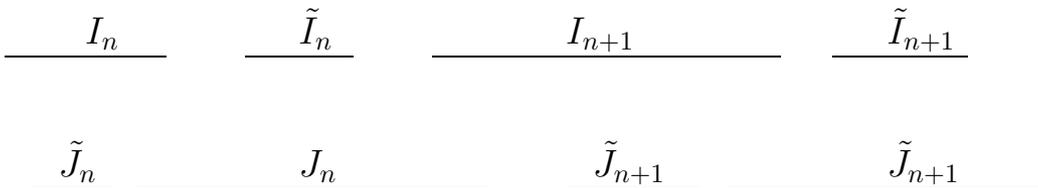
\begin{picture}(0, 0)

\put(-10,50){  \line(1, 0){60}}
\put(80,50){  \line(1, 0){40}}
\put(150,50){  \line(1, 0){130}}
\put(300,50){  \line(1,0){50}}

\put(20,55){ $ I_n$}\put(100,55){ $\tilde  I_n$}\put(200,55){ $ I_{n+1}$}\put(320,55){ $\tilde I_{n+1}$}
\put(80,50){  \line(1, 0){40}}

\put(10,5){ $ \tilde J_n$}\put(100,5){ $  J_n$}\put(210,5){ $ \tilde J_{n+1}$}\put(320,5){ $\tilde J_{n+1}$}
\put(80,50){  \line(1, 0){40}}

\put(40,0){  \line(1, 0){130}}
\put(0,0){  \line(1, 0){30}}
\put(200,0){  \line(1, 0){50}}
\put(260,0){  \line(1,0){120}}


\end{picture}

\vspace{2mm} 
Slightly more realistic picture:
\\ \\

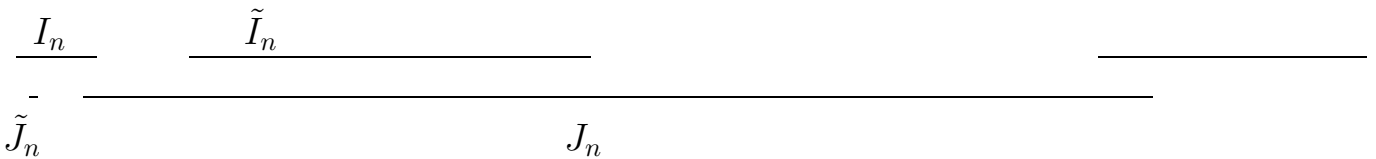
\begin{picture}(0, 0)

\put(0,0){  \line(1, 0){3}}
\put(-5,15){  \line(1, 0){30}}
\put(60,15){  \line(1, 0){150}}
\put(20,0){  \line(1, 0){400}}
\put(400,15){  \line(1, 0){100}}

\put(1,20){ $ I_n$}
\put(-10,-20){ $\tilde  J_n$}
\put(200,-20){ $  J_{n}$}
\put(80,20){ $\tilde I_{n}$}

\end{picture}
\\ \\ \\
Rank-one transformations $S$,$T$ via large  spacers could possess the following property (see  	arXiv:1106.4655  for rank-one  infinite transformations with large spacers).  Given functions $f$,$g$ with the  supports in fixed towers, respectively, for all large $k$ one has 
 the following zero-correlations:
$$(f, S^nf)=0, \  n\in I_k,$$ 
$$(g, T^ng)=0,  \ n\in J_k.$$
 As a consequence we get  for all large $n$
$$ \left(f\otimes g, S^n f\otimes T^n g\right)=0,$$
and, in fact, this is a case  for all $f,g$ from  a   dense family of vectors.
Thus, \it the product $S\otimes T$  has Lebesgue spectrum. \rm

We use time intervals $\tilde I_n, \tilde J_n$, $ \tilde I_n\subset J_n,$ $  \tilde J_n\subset I_n$,    to provide  the generic  behavior
of our transformations $S,T$.
So,  $S$,$T$ could  be  rigid, of  simple spectrum as well as their 
Gaussian and Poisson suspensions. We see that   constructions of 
\it  rigid 
Poisson   suspensions  $S_\ast$ and $T_\ast$, for which  
$S_\ast\otimes T_\ast$ has Lebesgue spectrum in $L_2^0\otimes L_2^0$, \rm
is reduced  to a  simple  rank-one exercise.  In addition, $S_\ast ,  T_\ast$ could be   of simple spectrum.  We get the same properties for associated Gaussian automorphisms $e^S,\ e^T$. Summarize and formulate

\bf THEOREM 1. \it  There are  rigid Poisson suspensions  $T_1$,$T_2$ with simple spectrum such that  $T_1\otimes T_2$ has  Lebesgue  spectrum
in the  orthocomplement to the sum of the coordinate spaces. There are rigid Gaussian transformations  $T_1'$,$T_2'$ with the same spectra as ones of 
$T_1$,$T_2$.
\rm 
\vspace{5mm} 
\\
{\Large \bf Let $S$ and $T$ act on a  probability space.}  
Below we  follow  author's draft which has been made  public  in 2009. 

 Given  vectors $h\in L_2^0$,  $\|h\|=1$,  
and  an  automorphism $R$  and $I\subset \N$,   denote  
$$Corr(R, h, I):=\sum_{n\in I}
| (R^nh, h)|.$$
\\
\bf LEMMA 1. \it  Given $\eps >0$ there are  generic maps $S$,$T$, vectors $f,g\in L_2^0(X,\mu)$,  $\|f\|=\|g\|=1$,  and integer intervals $I_n, J_n$,
$n=1,2, \dots ,$
such that  $$\bigcup_n  (I_n\cup J_n)=\N.  $$
$${ Corr}(S,  f, I_k) <  \frac \eps {2^{k+1}}, \ \ \
Corr(T, g,J_k) <  \frac \eps {2^{k+1}},$$
so 
$$Corr(S\otimes T,  f\otimes g, \N)<\eps.$$\rm
\\
\bf THEOREM 2. \it  There are  "generic" maps $S$,$T$ such that  $S\otimes T$ has  Lebesgue component in spectrum.
\rm

Theorem 2  follows obviously from Lemma 1.

Define
$\rho(T, \tilde T) =\mu ( supp (T^{-1}\tilde T).$
\\
\bf LEMMA 2. a)   \it Given $T$ with Lebesgue spectrum and  $\delta >0$, there exists  generic 
$\tilde T$ such that $\rho(T, \tilde T)<\delta$.

b)   \it Given generic $\tilde T$ and  $\delta >0$, there exists  an automorphism 
$ T$  with Lebesgue spectrum  such that $\rho(T, \tilde T)<\delta$.
\rm 

Lemma 2 is a consequence of the classic Rokhlin-Kakutani lemma.
\\
\bf LEMMA 3.  \it Let $T$  have  Lebesgue spectrum. Given $f\in L_2^0(X,\mu)$,  $\|f\|=1$, and  $\delta >0$, there exist $M$ and  
$f'\in L_2^0(X,\mu)$,  $\|f'\|=1$,  $\|f-f'\|<\delta$ such that 
$$(T^mf', f')=0$$ for all  $m>M$.\rm

Proof.  Consider an orthogonal base
$ h_{i,j},   \ \ Th_{i,j}=h_{i+1,j},$
find a finite linear combination $f'$,  $\|f-f'\|<\delta$, then set $M=1+max \{\ |i|\ :\
f'\  is \ not \ orthogonal  \ to \ h_{i,j}\}$.

 \rm  To build  transformations
let's apply the following "algorithm".
We consider sequences of automorphisms and functions
$T_k, \tilde T_k,  \ S_k, \tilde S_k,  \ f_k, g_k$, $\|f_k\|=\|g_k'\|=1$,
such that $S_k, T_k$ have Lebesgue spectrum, $\tilde T_k,  \tilde S_k$
are generic, 
 $$\rho(T_k, \tilde T_k)\to 0,  \ \ \rho(\tilde T_k, T_{k+1})\to 0,$$
$$\rho(S_k, \tilde S_k)\to 0,  \ \ \rho(\tilde S_k, S_{k+1})\to 0,$$
$$  f_k\to f, \ \ g_k\to g, \ \ k\to\infty$$
where the rate of the  convergences is chosen as fast as we want to provide 
$T_k\to T$, $S_k\to S$  with 
$$ Corr (T_k,f_k, J_k)=0,  \  Corr (S_k,f_k, I_k)=0,$$
$$ Corr (S,f, I_k)<\frac \eps {2^{k+1}},\ \
Corr (T,g, I_k)<\frac \eps {2^{k+1}}.$$
Summarizing the above, we obtain
$$Corr(S\otimes T,  f\otimes g, \N)<\eps.$$
Thus, $S\otimes T$ has the Lebesgue component in spectrum.

\bf Generic behavior \rm of $T$  is inherited  by our $\tilde T_k$.  
On the union $\bigcup \tilde J_n$   the behavior of powers $T^n$ 
is generic:  for instance,  we can  get all polynomial limits as $n$ runs
  within   $\bigcup \tilde J_n$.
In addition we provide   $T$ to be rank-one.  The same thing for $S$.
So,  $T$ and $S$  will be of simple spectrum, moreover  the same property hold   for  all their symmetric powers.  Hence,
Gaussian  "suspensions"  $e^{S}$ and $e^{T}$ are rigid, have simple spectrum.
  Their spectra are "very" singular,  but the product  $e^{S}\otimes e^{T}= e^{S\oplus T}$  has Lebesgue component in spectrum
(we do not prove that  in fact it  is  Lebesgue in $L_2^0\otimes L_2^0$).

{
\bf REMARK. \rm  
 A  similar effect  via different  methods is  presented by 
V. Bergelson, A. del Junco, M. Lemanczyk, J. Rosenblatt in arXiv:1103.0905.
There the authors assert  (Proposition 3.35)  \bf for arbitrary \rm 
 rigid, weakly mixing dynamical system $T$ the existence of a such $S$ that
in   the orthocomplement  of the sum of  multiplier spaces 
the product  $T\otimes S$  has absolutely continuous spectrum.}
\end{document}